\documentclass{article}

\usepackage[utf8]{inputenc}
\usepackage{authblk}
\usepackage{setspace}
\usepackage[margin=1.25in]{geometry}
\usepackage{graphicx}
\graphicspath{ {./figures/} }
\usepackage{subcaption}
\usepackage{amsmath}
\usepackage{lineno}
\usepackage{amsfonts}
\usepackage{amssymb, amsthm, amsmath, amsfonts}
\usepackage{wasysym}
\usepackage{mathrsfs}
\usepackage{hyperref}
\usepackage{graphicx}
\usepackage{lineno}
\usepackage[colorinlistoftodos]{todonotes}
\usepackage{listings}
\usepackage{cancel, enumerate}
\usepackage{rotating, environ}
\usepackage{caption}
\usepackage{subcaption}
\usepackage[inline]{enumitem}
\usepackage{dirtree}
\usepackage{xcolor}
\usepackage{tipa}
\usepackage{float}
\usepackage{changepage}
\usepackage{amsmath,amsthm,amssymb,amsfonts,mathrsfs,dsfont,ragged2e,tikz,tikz-cd,skull,color,enumitem,yfonts,graphicx,tkz-euclide,wasysym,fancyhdr,mathtools}

\newtheorem{defn}{Definition}
\newtheorem{prop}{Proposition}
\newtheorem{lemma}{Lemma}
\newtheorem{cor}{Corollary}
\newtheorem{ex}{Example}

\newtheorem{rmk}{Remark}

\newtheorem{thm}{Theorem}
\newtheorem*{ack}{Acknowledgements}

\renewcommand{\P}{\mathbb{P}}
\newcommand{\subseq}{\subseteq}

\newcommand{\Z}{\mathbb{Z}}

\newcommand{\Q}{\mathbb{Q}}

\newcommand{\R}{\mathbb{R}}

\newcommand{\Aut}{\text{Aut}}

\newcommand{\C}{\mathbb{C}}
\newcommand{\F}{\mathbb{F}}
\newcommand{\Hom}{\text{Hom}}

\newcommand{\id}{\text{id}}

\renewcommand{\P}{\mathbb{P}}

\newcommand{\wilde}{\widetilde}

\newcommand{\Gcal}{\mathcal{G}}
\newcommand{\Kcal}{\mathcal{K}}
\newcommand{\Lcal}{\mathcal{L}}

\renewcommand{\d}{\mathfrak{d}}

\newcommand{\Ecal}{\mathcal{E}}

\newcommand{\inv}{^{-1}}

\makeatletter
\newcommand{\tpitchfork}{%
  \vbox{
    \baselineskip\z@skip
    \lineskip-.52ex
    \lineskiplimit\maxdimen
    \m@th
    \ialign{##\crcr\hidewidth\smash{$-$}\hidewidth\crcr$\pitchfork$\crcr}
  }%
}
\makeatother
\newcommand{\PGL}{\text{PGL}}

\newcommand{\Pcal}{\mathcal{P}}
\newcommand{\Qcal}{\mathcal{Q}}

\newcommand{\Gr}{\text{Gr}}
\newcommand{\Dcal}{\mathcal{D}}

\newcommand{\wX}{\wilde{X}}
\newcommand{\wY}{\wilde{Y}}
\newcommand{\bL}{\mathbf{L}}

\usepackage[style=numeric,citestyle=numeric-comp,backend=biber,maxbibnames=9]{biblatex}
\title{Finite groupoids of configurations of lines in $\P^3_\C$}

\author{Jake Kettinger}

\affil{Department of Mathematics, Colorado State University}

\date{}

\begin{document}

\maketitle

\begin{abstract}
In this paper, we investigate groupoids coming from configurations of lines in three-dimensional space. Given a point and two skew lines in $\P^3_K$ over a field $K$, there exists a unique line containing the given point and meeting the two given lines. We use this construction to define a projection function from one line to another by using a skew line as an auxiliary. This way, we may create a groupoid whose objects are lines in a configuration, and whose morphisms are induced by these projection functions. We look at specific configurations for $K=\C$ that yield groupoids with finite automorphism groups.
\end{abstract}


\section{Introduction \& History}
While the study of configurations of points and lines is quite old indeed \cite{Reye1876, Kantor1879}, the study of groupoids induced by configurations of lines in $\P^3$ is quite recent. Although groupoids in general are fundamental to algebraic topology \cite{Brown}, the earliest mention of the form of groupoid studied in this paper comes from an unpublished 2008 article by Zakharevich \cite{ilya}. After that, the earliest discussion of groupoids induced by configurations of lines in three-dimensional space is found in \cite{politusgroupoid} and Ganger's Ph.D. thesis \cite{Ganger}.

The modern study of groupoids of configurations of lines grew out of the still quite recent study of geproci configurations of points and arrangements having the unexpected cone property. The geproci property was first explored back in 2020 by Pokora, Szemberg, and Szpond \cite{PSS}. Since then, a wealth of knowledge has been produced as Chiantini, Farnik, Favacchio, Harbourne, and Migliore joined in the study of geproci configurations \cite{politus1, politus2, politus3}. Since then, the geproci property has been studied in positive characteristic \cite{Kettinger, Ganger}, which has lead to a prominent connection to particular configurations of lines known as spreads. The study of spreads in the positive characteristic setting has enjoyed a long history \cite{BB, M} and has widespread applications coding theory \cite{Gorla}.

The main results of this paper are Theorems \ref{D4S3}, \ref{PS4}, and \ref{KS4}, found in Sections 3, 4, and 5, wherein we prove that the groupoids induced by the configurations (and subconfigurations) of lines spanning the $D_4$ configuration, Penrose configuration, and Klein configuration all have finite automorphism groups. We additionally determine which finite subgroup of $\PGL(2,\C)$ each group is isomorphic to.

\begin{ack}
We would like to thank Brian Harbourne and Allison Ganger for introducing us to this topic. We thank Chris Peterson for his guidance and insights while discussing this problem.
\end{ack}
\section{Background}
We will begin with some background on what a groupoid is and what is meant by the automorphism group of a groupoid.
\begin{defn}
    A \textbf{groupoid} is a category $\Gcal$ whose every morphism is invertible.
\end{defn}
Given a groupoid $\Gcal$ and and an object $A\in\Gcal$, the hom-set $\Hom_\Gcal(A,A)=\Aut_\Gcal(A)$ is a group. 

\begin{defn}
If $\Gcal$ is a groupoid and for all objects $A,B\in\Gcal$ we have $\Aut_\Gcal(A)\cong\Aut_\Gcal(B)$, then $\Gcal$ is a \textbf{regular} groupoid.
\end{defn}
\begin{rmk}
    For two objects $A,B$ in a groupoid $\Gcal$, $\Aut_\Gcal(A)\cong\Aut_\Gcal(B)$ if $\Hom_\Gcal(A,B)\neq\varnothing$. Thus a groupoid where every hom-set between two objects is nonempty must be regular. Such a groupoid is called \textbf{connected}.
\end{rmk}
\begin{defn}
In the case a groupoid $\Gcal$ is regular, we will denote the isomorphism class of all automorphism groups of the objects of $\Gcal$ as simply $\Aut_\Gcal$. We will refer to this group as \textbf{the group of the groupoid}.
\end{defn}
Now let $L_1,L_2,L_3\subseq\P^3_k$ be lines over some field $k$, where $L_1\cap L_2=L_2\cap L_3=\varnothing$. We define the isomorphism $$\pi(L_1,L_2,L_3):L_1\to L_3$$ as follows: Let $P\in L_1$ be a point. There is a unique plane $H$ that contains $P$ and $L_2$, since $L_1\cap L_2=\varnothing$. Then $\pi(L_1,L_2,L_3)(P):=L_3\cap H$, which is a point because $L_2\cap L_3=\varnothing$.

Note that it is only required that $L_1\cap L_2=L_3\cap L_2=\varnothing$. The map $\pi(L_1,L_2,L_3)$ is still well-defined and invertible in the case $L_1\cap L_3\neq\varnothing$. We will often refer to the line $L_2$ as the ``auxiliary line" of the map $\pi(L_1,L_2,L_3)$.

We can also refer to $\pi(L_1,L_2,L_3)$ with the indices of the lines involved, as follows: $\pi(1,2,3):=\pi(L_1,L_2,L_3)$. Note that since $\pi(i,j,k)=\pi(k,j,i)\inv$, we have $$\pi(k,j',i)\circ\pi(i,j,k)=(\pi(k,j,i)\circ\pi(i,j',k))\inv.$$

\begin{prop}
    Let $\Lcal$ be a set of lines in $\P^3$. Denote by $\Pi$ the set of all compositions of functions of the form $\pi(L_1,L_2,L_3)$, with $L_1,L_2,L_3\in\Lcal$ and $L_1\cap L_2=L_2\cap L_3=\varnothing$. Then $(\Lcal,\Pi)$ is a groupoid.
\end{prop}
Note that when $\Aut_{(\Lcal,\Pi)}$ is finite, the orbits of the groupoid are geproci \cite{politusgroupoid}.  This paper aims to examine configurations that are known to be geproci and determine whether they are orbits of a finite groupoid.
\begin{defn}
    A \textbf{parametrized line} in $\P^3$ is a pair $(L,f)$ of a line $L\subseq\P^3$ and a parametrization $f:\P^1\to L$.
\end{defn}
Given a line $L$, choosing a parametrization is equivalent to choosing a basis for $L$; one may define $L$ as the span of two points $P,Q\in\P^3$. If we write any point on $L$ in the form $aP+bQ$ for $(a,b)\in\P^1$, then the induced parametrization $f:\P^1\to L$ satisfies $f(a,b)=aP+bQ$, and vice versa.

Hence we can write the parametrized line $(L,f)$ as $\overline{PQ}$.

\begin{defn}
    Let $U=\overline{u_0u_1},V=\overline{v_0v_1},W=\overline{w_0w_1}$ be lines in $\P^3$. Then any point on $U$ can be written as $au_0+bu_1$ for $(a,b)\in\P^1$ and any point on $W$ can be written as $cw_0+dw_1$ for $(c,d)\in\P^1$. Then if $cw_0+dw_1=\pi(U,V,W)(au_0+bu_1)$, we must have $\overline{(au_0+bu_1)(cw_0+dw_1)}\cap V\neq\varnothing$.

    Therefore the wedge product $$(au_0+bu_1)\wedge v_0\wedge v_1\wedge (cw_0+dw_1)=0.$$ Expanding this out, we have $$ac(u_0\wedge v_0\wedge v_1\wedge w_0)+ad(u_0\wedge v_0\wedge v_1\wedge w_1)+bc(u_1\wedge v_0\wedge v_1\wedge w_0)+bd(u_1\wedge v_0\wedge v_1\wedge w_1)=0.$$ We can write this as a matrix formula $$\begin{pmatrix}
        a&b\\
    \end{pmatrix}\begin{pmatrix}
        u_0\wedge v_0\wedge v_1\wedge w_0&u_0\wedge v_0\wedge v_1\wedge w_1\\
        u_1\wedge v_0\wedge v_1\wedge w_0&u_1\wedge v_0\wedge v_1\wedge w_1\\
    \end{pmatrix}\begin{pmatrix}
        c\\d\\
    \end{pmatrix}=0,$$ or rewrite this equivalently: $$\begin{pmatrix}
        -u_0\wedge v_0\wedge v_1\wedge w_1&-u_1\wedge v_0\wedge v_1\wedge w_1\\
        u_0\wedge v_0\wedge v_1\wedge w_0&u_1\wedge v_0\wedge v_1\wedge w_0\\
    \end{pmatrix}\begin{pmatrix}
        a\\b\\
    \end{pmatrix}=\begin{pmatrix}
        c\\d\\
    \end{pmatrix}.$$

    Note that for points $P,Q,R,S\in\P^3$, the wedge product $P\wedge Q\wedge R\wedge S$ is simply the scalar $\det\begin{pmatrix}
        P^{\mathsf{T}}&Q^{\mathsf{T}}&R^{\mathsf{T}}&S^{\mathsf{T}}\\
    \end{pmatrix}$. 
     Furthermore, $$\begin{pmatrix}
        w_0&w_1\\
    \end{pmatrix}\begin{pmatrix}
        -u_0\wedge v_0\wedge v_1\wedge w_1&-u_1\wedge v_0\wedge v_1\wedge w_1\\
        u_0\wedge v_0\wedge v_1\wedge w_0&u_1\wedge v_0\wedge v_1\wedge w_0\\
    \end{pmatrix}\begin{pmatrix}
        a\\b\\
    \end{pmatrix}=\pi(U,V,W)(au_0+bu_1),$$ so we may represent $\pi(U,V,W)$ with this $2\times 2$ matrix in $\PGL(2,\C)$, given bases for $U$, $V$, and $W$.
    \end{defn}
The matrix in $\PGL(2,\C)$ can also be defined as the unique map $\P^1\to\P^1$ causing the following diagram to commute, given parametrizations $f:\P^1\to U$ and $h:\P^1\to W$. 
\begin{center}
    \begin{tikzcd}
        \P^1\arrow[rr]\arrow[dd, "f"]&&\P^1\arrow[dd, "h"]\\
        &&\\
        U\arrow[rr, "\pi(U{,}V{,}W)"]&&W\\
    \end{tikzcd}
\end{center}
The following command is Macaulay2 code for computing the projection map in $\PGL(2,\C)$ of a parametrized triple \texttt{(i,j,k)} where \texttt{i}, \texttt{j}, and \texttt{k} are lists of two lists of four numbers, e.g. of the form \texttt{i=\{\{1,2,3,4\},\{5,6,7,8\}\}}.\\

\texttt{proj=(i,j,k)->matrix\{\{-det(matrix\{i\_0,j\_0,j\_1,k\_1\}),-det(matrix\{i\_1,j\_0,j\_1,k\_1\})\},\\\{det(matrix\{i\_0,j\_0,j\_1,k\_0\}),det(matrix\{i\_1,j\_0,j\_1,k\_0\})\}\}}
\begin{ex}
    Consider the set $\Lcal$ of the four parametrized lines in $\P^3$: $a=\overline{(1,0,0,0)(0,1,0,0)}$, $b=\overline{(0,0,1,0)(0,0,0,1)}$, $c=\overline{(1,0,1,0)(0,1,0,1)}$, and $d=\overline{(-1,0,1,0)(0,-1,0,1)}$. Note that each line in $\Lcal$ is skew to every other line in $\Lcal$. In this case, $\Aut_{(\Lcal,\Pi)}\cong \id$, as can be verified computationally. Without the use of computational assistance, one may conclude that the group of the groupoid is trivial because the four lines lie on a quadric surface; therefore, any point $P$ on each line $L\in\Lcal$ is on a transversal line and is thus fixed by $\Aut_{(\Lcal,\Pi)}(L)$ \cite{Ganger}.
    
%

\end{ex}
\begin{ex}
    Consider the set $\Lcal$ of the four parametrized lines in $\P^3$: $a=\overline{(1,0,0,0)(0,1,0,0)}$, $b=\overline{(0,0,1,0)(0,0,0,1)}$, $c=\overline{(1,0,1,0)(0,1,0,1)}$, and $d=\overline{(x_0,y_0,z_0,w_0)(x_1,y_1,z_1,w_1)}$. Then according to Macaulay2, $$\pi(c,d,a)\circ\pi(a,b,c)=\begin{pmatrix}
        w_0x_1-w_0z_1-x_0w_1+z_0w_1&-z_0x_1+x_0z_1\\
        w_0y_1-y_0w_1&-z_0y_1+y_0z_1-w_0z_1+z_0w_1\\
    \end{pmatrix}.$$ This matrix will have finite order in $\PGL(2,\C)$ if the quotient of its eigenvalues is a root of unity.
\end{ex}

\section{The $D_4$ Configuration}
The $D_4$ configuration (also known as the Reye configuration) is a set of 12 points in $\P^3$ with a set of 16 lines, each containing three of the points. It is thus a $(12_4,16_3)$-configuration, following the $(p_q,n_k)$ notation for configurations \cite{G}. Under this notation, $p$ is the number of objects of type $A$ (in this case, points), $n$ is the number of objects of type $B$ (in this case, lines), and $q$ (resp. $k$) is the number of objects of type $B$ (resp. $A$) coinciding with each object of type $A$ (resp. $B$). The $D_4$ configuration is furthermore known to be a $(3,4)$-geproci half-grid \cite{politus1}.

    Let $\Lcal$ be the 16 lines of the $D_4$ configuration. For example, we can set $L_i\in\Lcal$ to be the rowspan of the $2\times 4$ matrix $\bL_i$ from the following set of matrices.   \begin{align*}
    \bL_0&=\begin{pmatrix}
        0&0&0&1\\
        0&0&1&0
    \end{pmatrix}&
    \bL_1&=\begin{pmatrix}
        0&0&0&1\\
        0&1&0&0
    \end{pmatrix}&
    \bL_2&=\begin{pmatrix}
        0&0&0&1\\
        1&0&0&0\\
    \end{pmatrix}&
    \bL_3&=\begin{pmatrix}
        0&0&0&1\\
        1&1&1&2
    \end{pmatrix}\\
    \bL_4&=\begin{pmatrix}
        1&1&0&1\\
        0&0&1&0
    \end{pmatrix}&
    \bL_5&=\begin{pmatrix}
        1&1&0&1\\
        0&1&0&0
    \end{pmatrix}&
    \bL_6&=\begin{pmatrix}
        1&1&0&1\\
        1&0&0&0
    \end{pmatrix}&
    \bL_7&=\begin{pmatrix}
        1&1&0&1\\
        1&1&1&2\\
    \end{pmatrix}\\
    \bL_8&=\begin{pmatrix}
        1&0&1&1\\
        0&0&1&0\\
    \end{pmatrix}&
    \bL_9&=\begin{pmatrix}
        1&0&1&1\\
        0&1&0&0
    \end{pmatrix}&
        \bL_{10}&=\begin{pmatrix}
        1&0&1&1\\
        1&0&0&0
    \end{pmatrix}&
    \bL_{11}&=\begin{pmatrix}
        1&1&0&1\\
        1&1&1&2
    \end{pmatrix}\\
    \bL_{12}&=\begin{pmatrix}
        0&1&1&1\\
        0&0&1&0
    \end{pmatrix}&
    \bL_{13}&=\begin{pmatrix}
        0&1&1&1\\
        0&1&0&0\\
    \end{pmatrix}&
    \bL_{14}&=\begin{pmatrix}
        0&1&1&1\\
        1&0&0&0\\
    \end{pmatrix}&
    \bL_{15}&=\begin{pmatrix}
        0&1&1&1\\
        1&1&1&2
    \end{pmatrix}
\end{align*}
   This configuration stands out from the remaining configurations we will explore in this paper, as it is the only one where the lines are not all mutually-skew. It nevertheless induces a connected (and thus regular) groupoid. Let $L_i,L_k\in\Lcal$. Then the invertible map $\pi(i,j,k)$ is defined for exactly two $j\in\{0,\dots,15\}\setminus\{i,k\}$, regardless of whether $L_i\cap L_k=\varnothing$ or $L_i\cap L_k\neq\varnothing$. This is because, given any pair of lines $L_i,L_k\in \Lcal$, there are exactly two lines in $\Lcal$ that are skew with both $L_i$ and $L_k$.

\begin{thm}\label{D4S3}
    The group of the groupoid of the $D_4$ is isomorphic to $S_3$.
\end{thm}
\begin{proof}
    First note that we can use the letters $A$, $B$, and $C$ to label the vertices of $D_4$ with elements from the set $\{A,B,C\}\oplus (\Z/2\Z)^2$ such that there is a bijective correspondence between lines of the $D_4$ and triples of the form $$\{A+g,B+g',C+g'':g+g'+g''=0\}.$$

    Figure \ref{D4figure} below shows an example, using $(\Z/2\Z)^2=\langle z_1,z_2\rangle$.
    \begin{figure}[H]
    \begin{center}
        \includegraphics[scale=0.5]{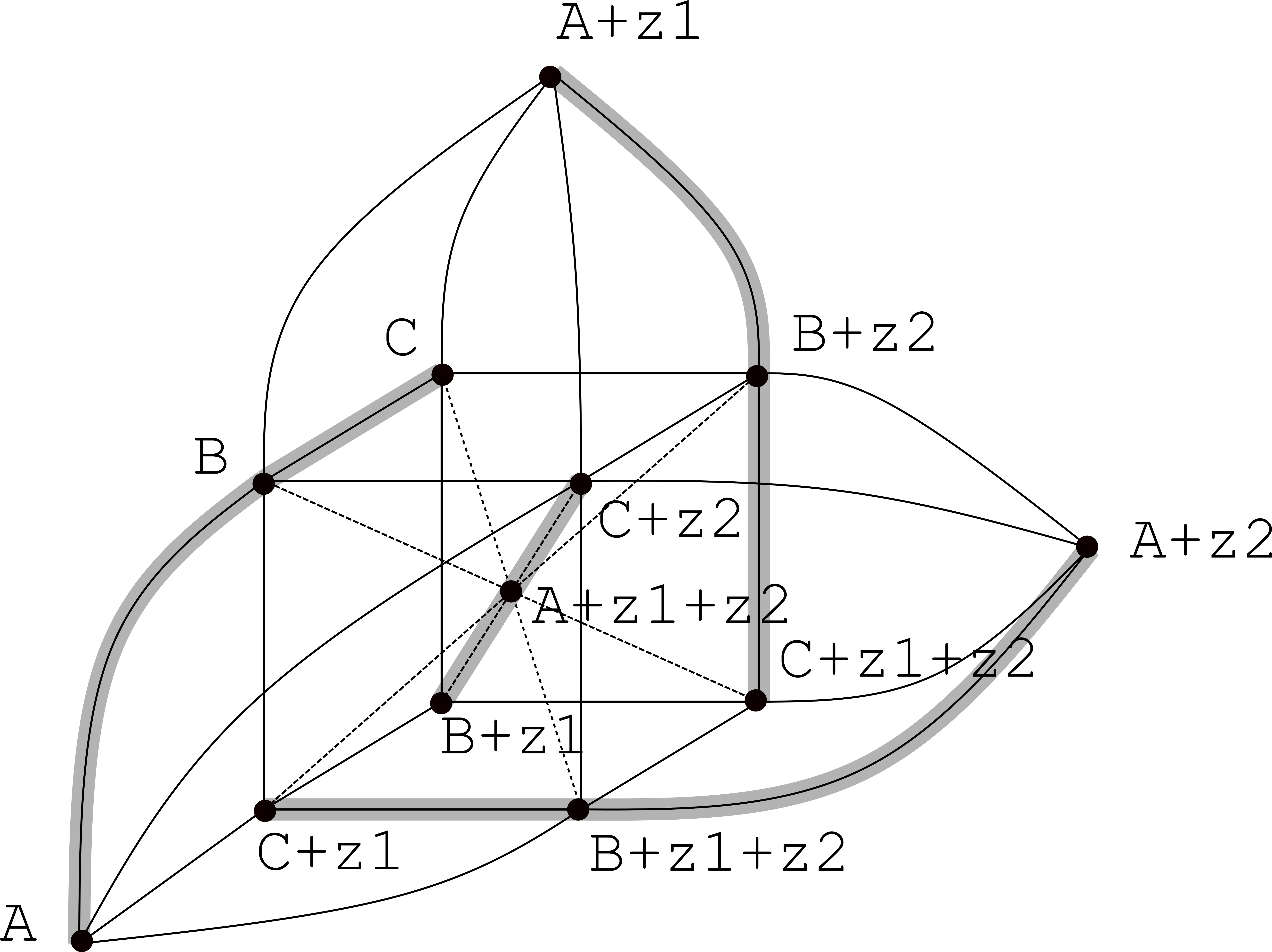}
        \end{center}
        \caption{A labeling of the $D_4$ configuration}
        \label{D4figure}
    \end{figure}
    The first step is to show that for every quadruple-point $Q$ of the configuration, and every $\pi\in\Pi$, that $\pi(Q)$ is another quadruple-point.

    Using the labeling above, we can define the functions $\pi(L_i,L_j,L_k)$ in terms of the triples. We can write 
    \begin{center}
    \begin{tikzcd}

    \{A+f_A,B+f_B,C+f_C\}\arrow[rrr,"\{A+g_A{,}B+g_B{,}C+g_C\}"]&&&\{A+h_A,B+h_B,C+h_C\}\\
    \end{tikzcd}
    \end{center}
    and determine where $A+f_A$, $B+f_B$, and $C+f_C$ must go. There are two cases: the case where $$\{A+f_A,B+f_B,C+f_C\}\cap \{A+h_A,B+h_B,C+h_C\}=\varnothing$$ and the case $$\{A+f_A,B+f_B,C+f_C\}\cap \{A+h_A,B+h_B,C+h_C\}\neq\varnothing.$$

    \underline{Case 1}: Since it is required that the auxiliary line not intersect the domain and codomain, we know that $f_A,g_A,h_A$ are distinct, $f_B,g_B,h_B$ are distinct, and $f_C,g_C,h_C$ are distinct. Then to find the image of $A+f_A$ under $\pi$, we want to show that either $$f_A+g_B+h_C=0$$ or $$f_A+g_C+h_B=0.$$ We can rewrite these as $$f_A+g_B+h_A+h_B=0$$ and $$f_A+g_A+g_B+h_B=0.$$ And so we would like to show that either $$f_A+g_B+h_B=h_A$$ or $$f_A+g_B+h_B=g_A.$$ Since $g_B\neq h_B$, we have $g_B+h_B\neq 0$, and so we know that $f_A+g_B+h_B\neq f_A$. Since there are only four elements of $(\Z/2\Z)^2$, that leaves the possibilities that $f_A+g_B+h_B$ is $h_A$, $g_A$, or $f_A+g_A+h_A$. We want to prove $f_A+g_B+h_B\neq f_A+g_A+h_A$.

    If $f_A+g_B+h_B=f_A+g_A+h_A$, then $g_B+h_B=g_A+h_A$. But then $g_B+g_A=h_B+h_A$ and so $g_C=h_C$, which is a contradiction. Therefore $f_A+g_B+h_B\neq f_A+g_A+h_A$. And so it can only be that $f_A+g_B+h_B=h_A$ or $f_A+g_B+h_B=g_A$. Then suppose without loss of generality that $f_A+g_B+h_B=h_A$. Then $f_A+g_B+h_C=0$ and so $\{A+f_A,B+g_B,C+h_C\}$ is a line containing $A+f_A$ and intersects the auxiliary line at $B+g_B$ and the codomain at $C+h_C$. Thus $\pi(A+f_A)=C+h_C$, and we are done with Case 1.

    \underline{Case 2}: Now suppose $\{A+f_A,B+f_B,C+f_C\}\cap \{A+h_A,B+h_B,C+h_C\}\neq\varnothing$. Without loss of generality, say $f_A=h_A$. Then in order for the two lines to be distinct, it is required that $f_B\neq h_B$ and $f_C\neq h_C$. Then the plane $H$ that contains the point $A+f_A$ and the line $\{A+g_A,B+g_B,C+g_C\}$ meets the line $\{A+h_A,B+h_B,C+h_C\}$ at the point $A+h_A=A+f_A$, so $\pi(A+f_A)=A+f_A$.

    Then we claim that $f_B+g_A+h_C=0$ and $f_C+g_A+h_B=0$, and so $\pi(B+f_B)=C+h_C$ and $\pi(C+f_C)=B+h_B$.

    We can rewrite the first equation as $f_A+f_C+g_A+h_A+h_B=0$, and since $f_A+h_A=0$, this is $f_C+g_A+h_B=0$, so the prior two equations are equivalent.

    We want to show that $f_B+h_C=g_A$. First note that $f_B+h_C\neq f_A$ because if so, then $f_A+f_B+h_A+h_B=f_B+h_B=0$, which is a contradiction.

    We want to prove the following three facts:
    \begin{enumerate}
        \item\label{fact1} $f_B+h_C+g_B\neq f_C$. Suppose $f_B+h_C+g_B= f_C$. Then $f_A+h_C+g_B=0$, and so $h_A+h_C+g_B=0$, and so $h_B+g_B=0$ and so $g_B=h_B$, a contradiction.

        \item\label{fact2} $f_B+h_C+g_B\neq h_C$. This is true because $f_B+g_B\neq 0$ since $f_B\neq g_B$.

        \item\label{fact3} $[[f_B+h_C+g_B]+[f_C]+[h_C]]+g_B=f_A$. This works out by cancellation and the equation $f_B+f_C=f_A$.
    \end{enumerate}
    What facts \ref{fact1} and \ref{fact2} show is that $f_B+h_C+g_B$ is a candidate for $g_C$, since it is a requirement that $g_C$ is distinct from both $f_C$ and $h_C$. The only two elements of $(\Z/2\Z)^2$ that are distinct from $f_C$ and $h_C$ are $f_B+h_C+g_B$ and $[f_B+h_C+g_B]+[f_C]+[h_C]$. This is because $[f_B+h_C+g_B]+[f_C]\neq 0$, $[f_B+h_C+g_B]+[h_C]\neq 0$, and $[f_C]+[h_C]\neq 0$.
    
    Fact \ref{fact3} shows that $g_C\neq [f_B+h_C+g_B]+[f_C]+[h_C]$ because then $g_A$ would be $g_C+g_B=[[f_B+h_C+g_B]+[f_C]+[h_C]]+g_B=f_A$, which is a contradiction because $f_A$ and $g_A$ must be distinct.

    And so we must have $g_C=f_B+h_C+g_B$. Then we have $f_B+h_C=g_B+g_C=g_A$, as required. Similarly, we have $f_C+h_B=g_A$. This means that $\{A+g_A,B+f_B,C+h_C\}$ is a line intersecting the domain, codomain, and auxiliary line and so $\pi(B+f_B)=C+h_C$. Similarly, $\pi(C+f_C)=B+f_B$.

    Therefore under every map in $\Pi$, every quadruple point maps to another quadruple point. Therefore the automorphism groups $\Aut_\Gcal(L)$ permute the three quadruple points of $L$. Since an automorphism of $\P^1$ is uniquely determined by where it sends three points, we can say that $\Aut_\Gcal(L)\leq S_3$. Since we already found automorphisms of order 2 and 3 above, we can then say $\Aut_\Gcal(L)\cong S_3$.
\end{proof}
\begin{rmk}
Only 6 lines of $\Lcal$ are required to induce the group $S_3$ from the groupoid: take four lines that form a square (like $L_4$, $L_8$, $L_{14}$ and $L_{15}$ in the above example) and a line that is skew to all four lines of the square (like line $L_1$). Then we can make a cycle of size 2. To those five lines, add a line that is skew to $L_1$, and two parallel lines of the square (for example, $L_{15}$ and $L_4$, which would yield line $L_{10}$). Adding $L_{10}$ to the configuration yields a mutually-skew subconfiguration ($L_1$, $L_4$, $L_{10}$ and $L_{15}$) which yields a cycle of order 3. Thus we can form all of $S_3$ with $\{L_1,L_4,L_8,L_{10},L_{14},L_{15}\}$.
\end{rmk}

\section{The Penrose Configuration}
\subsection{The Groupoid of the half-Penrose Spread}
Let $t\in\C$ satisfy $t^2+t+1=0$. Then the Penrose configuration is the following set of 40 points \cite{VlasovP}.

\begin{table}[htb]
\[
\begin{array}{|l|l|l|l|}
\hline
 (1,\ 0,\ 0,\ 0)& (0,\ 1,\ 0,\ 0)&
\ ( 0,\ 0,\ 1,\ 0)& (0,\ 0,\ 0,\ 1)\\ \hline
 (0,\ 1,-1,\ 1)& (1,\ 0,-1,-1)&
(1,-1,\ 0,\ 1)& (1,\ 1,\ 1,\ 0)\\ \hline
(0,\ 1,-t,\ t^2)&
(1,\ 0,-t,-t^2)&
(1,-t,\ 0,\ t^2)&
(1,\ t,\ t^2,\ 0)\\ \hline (0,\ 1,-t^2,\ t)&
(1,\ 0,-t^2,-t)&
(1,-t^2,\ 0,\ t)&
(1,\ t^2,\ t,\ 0)\\ \hline
 (0,\ 1,-t,\ 1)&
(1,\ 0,-1,-t)&
(1,-t^2,\ 0,\ t^2)&
(1,\ t,\ 1,\ 0)\\ \hline
 (0,\ 1,-t^2,\ t^2)&
(1,\ 0,-t,-1)&
(1,-1,\ 0,\ t)&
(1,\ t^2,\ t^2,\ 0)\\ \hline
 (0,\ 1,-1,\ t)&
(1,\ 0,-t^2,-t^2)&
(1,-t,\ 0,\ 1)&
(1,\ 1,\ t,\ 0)\\ \hline
 (0,\ 1,-t^2,\ 1)&
(1,\ 0,-1,-t^2)&
(1,-t,\ 0,\ t)&
(1,\ t^2,\ 1,\ 0)\\ \hline
 (0,\ 1,-1,\ t^2)&
(1,\ 0,-t,-t)&
(1,-t^2,\ 0,\ 1)&
(1,\ 1,\ t^2,\ 0)\\ \hline
 (0,\ 1,-t,\ t)&
(1,\ 0,-t^2,-1)&
(1,-1,\ 0,\ t^2)&
(1,\ t,\ t,\ 0) \\ \hline
 \end{array}
\]
\label{tab40rays}
\end{table}

The Penrose configuration $\Pcal$ is known to be a $(5,8)$-geproci non-half-grid \cite{politus1}. Furthermore, $\Pcal$ is a $(40_9,90_4)$-configuration of points and lines.

  We may define a set $\Lcal$ of ten mutually-skew lines $L_0,\dots,L_9$ defined as the rowspans of the respective matrices $\mathbf{L}_i$, given as follows.

\begin{align*}
    \bL_0&=\begin{pmatrix}
        1&0&-t^2&-t^2\\
        1&0&-t&-t
    \end{pmatrix}&
    \bL_1&=\begin{pmatrix}
        0&1&-t^2&t^2\\
        0&1&-t&t
    \end{pmatrix}\\
    \bL_2&=\begin{pmatrix}
        1&0&-t^2&-t\\
        0&1&-t&t^2\\
    \end{pmatrix}&
    \bL_3&=\begin{pmatrix}
        0&1&-t^2&t\\
        1&0&-t&-t^2
    \end{pmatrix}\\
    \bL_4&=\begin{pmatrix}
        0&1&-t&1\\
        1&0&-1&-t
    \end{pmatrix}&
    \bL_5&=\begin{pmatrix}
        1&0&-t&-1\\
        0&1&-1&t
    \end{pmatrix}\\
    \bL_6&=\begin{pmatrix}
        1&0&-1&-t^2\\
        0&1&-t^2&1
    \end{pmatrix}&
    \bL_7&=\begin{pmatrix}
        0&1&-1&t^2\\
        1&0&-t^2&-1\\
    \end{pmatrix}\\
    \bL_8&=\begin{pmatrix}
        1&1&t&0\\
        1&1&t^2&0\\
    \end{pmatrix}&
    \bL_9&=\begin{pmatrix}
        1&-1&0&t\\
        1&-1&0&t^2
    \end{pmatrix}\\
\end{align*}

Each line $L_i$ contains four of the points of $\Pcal$. Furthermore, for $i\neq8,9$, the basis chosen for $L_i$ yields that the four points of $\Pcal\cap L_i$ are given by the parameters $$\begin{pmatrix}
    1&0
\end{pmatrix},\,\begin{pmatrix}
    0&1
\end{pmatrix},\,\begin{pmatrix}
    t&1
\end{pmatrix},\,\begin{pmatrix}
    -t^2&1
\end{pmatrix}.$$ For example, the four points of $L_0\cap\Pcal$ are \begin{align*}
    (1,0,-t^2,-t^2)&=1\cdot(1,0,-t^2,-t^2)+0\cdot(1,0,-t,-t),\\
    (1,0,-t,-t)&=0\cdot(1,0,-t^2,-t^2)+1\cdot(1,0,-t,-t),\\
    (1,0,-1,-1)&=t\cdot(1,0,-t^2,-t^2)+1\cdot(1,0,-t,-t),\\
    (1,0,0,0)&=-t^2\cdot(1,0,-t^2,-t^2)+1\cdot(1,0,-t,-t).\\
\end{align*}
For $L_8$ and $L_9$, the Penrose points are instead given by the parameters $$\begin{pmatrix}
    1&0
\end{pmatrix},\,\begin{pmatrix}
    0&1
\end{pmatrix},\,\begin{pmatrix}
    -1&1
\end{pmatrix},\,\begin{pmatrix}
    t^2&1
\end{pmatrix}.$$

\begin{lemma}
    Let $X=\left\{\begin{pmatrix}
        1\\0
    \end{pmatrix},\begin{pmatrix}
        0\\1
    \end{pmatrix},\begin{pmatrix}
        t\\1
    \end{pmatrix},\begin{pmatrix}
        -t^2\\1
    \end{pmatrix}\right\}\subseq\P^1_\C$, and $\Aut_X(\P^1_\C)=\{f\in\Aut(\P^1_\C):f(X)= X\}$. Then $\Aut_X(\P^1)\cong A_4$. In other words, the four points of $X$ have tetrahedral symmetry.
\end{lemma}
\begin{proof}
    Since $\Aut_X(\P^1)$ acts on the finite set of points $X$ with $\# X\geq3$, we must have $\Aut_X(\P^1)$ be a finite group. This is because three conditions uniquely determine an automorphism of $\P^1$, so no two distinct automorphisms of $\P^1$ will have the same permutation action on $X$. Thus $\Aut_X(\P^1)\to S_4$ mapping $f$ to its permutation on the points of $X$ is an embedding, so $\Aut_X(\P^1)$ is isomorphic to a subgroup of $S_4$.
    
    The only finite subgroups of $\Aut(\P^1_\C)=\PGL(2,\C)$ are (up to isomorphism): $C_n$, the cyclic group of size $n$, for any $n\in\Z$; $D_{2n}$, the dihedral group of size $2n$, for any $2n\in2\Z$; $A_4$; $S_4$; and $A_5$ \cite{Klein}.

    Furthermore, $\Aut_X(\P^1)$ must be isomorphic to a subgroup of $S_4$, so $\Aut_X(\P^1)\not\cong A_5$.

    We claim that \begin{align*}
    \Aut_X(\P^1)=G:=&\left\{\begin{pmatrix}
        1&0\\0&1
    \end{pmatrix},\begin{pmatrix}
        0&-1\\1&0
    \end{pmatrix},\begin{pmatrix}
        t&1\\1&-t
    \end{pmatrix},\begin{pmatrix}
        -1&t\\t&1
    \end{pmatrix},\begin{pmatrix}
        0&-1\\t&-t^2
    \end{pmatrix},\begin{pmatrix}
        t&-t^2\\1&0
    \end{pmatrix}, \right.\\&\left.\begin{pmatrix}
        0&-t\\1&t^2
    \end{pmatrix},\begin{pmatrix}
        t^2&t\\-1&0
    \end{pmatrix},\begin{pmatrix}
        t&0\\1&t^2
    \end{pmatrix},\begin{pmatrix}
        t^2&0\\-1&t
    \end{pmatrix},\begin{pmatrix}
        t^2&t\\0&1
    \end{pmatrix},\begin{pmatrix}
        t&-t^2\\0&1
    \end{pmatrix}\right\}.
    \end{align*}

It is straightforward to verify computationally that every matrix in the group $G\leq\PGL(2,\C)$ permutes the points of $X\subseq\P^1$, so $G\leq\Aut_X(\P^1)$. Specifically, the action of $G$ on $X$ induces \textit{even} permutations of the points of $X$, giving us that $G\cong A_4$. Thus either $\Aut_X(\P^1)\cong A_4$ or $S_4$, since no cyclic group $C_n$ or dihedral group $D_{2n}$ contains a subgroup isomorphic to $A_4$.

We will show that $\Aut_X(\P^1)\not\cong S_4$ by demonstrating that there is no automorphism in $f\in\Aut(\P^1)$ performing the following odd permutation on $X$: \begin{align*}
    f\begin{pmatrix}
        1\\0
    \end{pmatrix}&=\begin{pmatrix}
        0\\1
    \end{pmatrix},&f\begin{pmatrix}
        0\\1
    \end{pmatrix}&=\begin{pmatrix}
        1\\0
    \end{pmatrix},\\
    f\begin{pmatrix}
        t\\1
    \end{pmatrix}&=\begin{pmatrix}
        t\\1
    \end{pmatrix},&f\begin{pmatrix}
        -t^2\\1
    \end{pmatrix}&=\begin{pmatrix}
        -t^2\\1
    \end{pmatrix}.\\
\end{align*}
The only matrix in $\PGL(2,\C)$ that satisfies the first three criteria is $\begin{pmatrix}
    0&t^2\\1&0
\end{pmatrix}$, which fails to satisfy the fourth criterion because $\begin{pmatrix}
    0&t^2\\1&0
\end{pmatrix}\begin{pmatrix}
    -t^2\\1
\end{pmatrix}=\begin{pmatrix}
    -1\\1
\end{pmatrix}$.
Therefore $\Aut_X(\P^1)\not\cong S_4$, since there exist odd permutations of $X$ that cannot be expressed as automorphisms of $\P^1$.

Note that three conditions uniquely determine an automorphism of $\P^1$, so there are no two distinct automorphisms $f\neq g\in\Aut(\P^1)$ such that $f$ and $g$ have the same action on $X$. Thus $\Aut_X(\P^1)=G\cong A_4$. 
\end{proof}
\begin{lemma}
    Let $Y=\left\{\begin{pmatrix}
        1\\0
    \end{pmatrix},\begin{pmatrix}
        0\\1
    \end{pmatrix},\begin{pmatrix}
        -1\\1
    \end{pmatrix},\begin{pmatrix}
        t^2\\1
    \end{pmatrix}\right\}\subseq\P^1_\C$. Then $\Aut_Y(\P^1)\cong A_4$.
\end{lemma}
\begin{proof}
    The proof is exactly the same as in Lemma 1; the only difference is \begin{align*}
    \Aut_Y(\P^1)=&\left\{\begin{pmatrix}
        1&0\\0&1
    \end{pmatrix},\begin{pmatrix}
        0&-1\\t&0
    \end{pmatrix},\begin{pmatrix}
        -1&t^2\\1&1
    \end{pmatrix},\begin{pmatrix}
        1&1\\t&-1
    \end{pmatrix},\begin{pmatrix}
        0&-1\\1&1
    \end{pmatrix},\begin{pmatrix}
        1&1\\-1&0
    \end{pmatrix}, \right.\\&\left.\begin{pmatrix}
        0&t\\-1&t^2
    \end{pmatrix},\begin{pmatrix}
        t^2&-t\\1&0
    \end{pmatrix},\begin{pmatrix}
        t&0\\-t&1
    \end{pmatrix},\begin{pmatrix}
        1&0\\t&t
    \end{pmatrix},\begin{pmatrix}
        t&-1\\0&1
    \end{pmatrix},\begin{pmatrix}
        1&1\\0&t
    \end{pmatrix}\right\}.
    \end{align*}
\end{proof}
\begin{cor}
    Let $L\in\Lcal$. Then $$\Aut_{L\cap\Pcal}(L)\cong A_4.$$
\end{cor}
\begin{proof}
    Follows from Lemmas 1 and 2.
\end{proof}
\begin{thm}\label{HPA4}
    Let $\Ecal=\{L_0,L_2,L_4,L_6,L_8\}$ and let $\Pi^*$ be the composition-closure of the set of morphisms of the form $\pi(L,L',L'')$, $L,L',L''\in\Ecal$ distinct. Then the groupoid $(\Ecal,\Pi^*)$ has the automorphism group $$\Aut_{(\Ecal,\Pi^*)}\cong A_4.$$
\end{thm}
Note that the 20 points of $\Pcal$ contained on the lines of $\Ecal$ are (isomorphic to) the half-Penrose configuration explored in Section 5.3 of \cite{politus1}. It and its complement in $\Pcal$ are $(4,5)$-geproci half-grids.
\begin{proof}
We begin by showing that for any $P\in\Pcal^*:=\bigcup_{L\in\Ecal}(\Pcal\cap L)$ and $\pi\in\Pi^*$, $$\pi(P)\in\Pcal^*.$$ First let us make use of the following definition.
\begin{defn}
Let $\Lcal$ be a configuration of lines in $\P^3$, and let $(\Lcal,\Pi)$ be the groupoid induced by the composition-closure of the projection maps between the lines of $\Lcal$. A morphism $\pi\in \Pi$ is \textbf{simple} if $\pi$ is of the form $\pi(L,L',L'')$ for distinct $L,L',L''\in\Lcal$.
\end{defn}
    The following Macaulay2 code shows that for any simple morphism $\pi(L,L',L'')\in\Pi^*$, $\pi(L,L',L'')(L\cap\Pcal)=L''\cap\Pcal$ \cite{M2}. Recall that lines $L_0$, $L_4$, $L_6$ and $L_8$ have Penrose points corresponding to the parameters $X$, and $L_2$ has its Penrose points at the parameters $Y$. In our Macaulay2 code, we define $E_i:=L_{2i}$.
    
    \texttt{i1: K=QQ[t]/ideal(t\^{}2+t+1)}\\

    \texttt{i2: E=\{\{\{1,0,t+1,t+1\},\{1,0,-t,-t\}\},\{\{1,0,t+1,-t\},\{0,1,-t,t\^{}2\}\},\\\{\{0,1,-t,1\},\{1,0,-1,-t\}\},\{\{1,0,-1,-t\^{}2\},\{0,1,-t\^{}2,1\}\},\\\{\{1,1,t,0\},\{1,1,t\^{}2,0\}\}\};}\\

\texttt{i3: proj=(i,j,k)->matrix\{\{-det(matrix\{i\_0,j\_0,j\_1,k\_1\}),\\-det(matrix\{i\_1,j\_0,j\_1,k\_1\})\},\{det(matrix\{i\_0,j\_0,j\_1,k\_0\}),\\det(matrix\{i\_1,j\_0,j\_1,k\_0\})\}\};}\\

\texttt{i4: X=\{matrix\{\{1\},\{0\}\},matrix\{\{0\},\{1\}\},\\matrix\{\{t\},\{1\}\},matrix\{\{-t\^{}2\},\{1\}\}\};}\\

\texttt{i5: Y=\{matrix\{\{1\},\{0\}\},matrix\{\{0\},\{1\}\},\\matrix\{\{-1\},\{1\}\},matrix\{\{t\^{}2\},\{1\}\}\};}\\

\texttt{i6:  Xsimples=select(flatten apply(\{0,1,2,3\},i->apply(\{0,1,2,3,4\},\\j->select(apply(\{0,1,2,3\},k->proj(E\_i,E\_j,E\_k)),m->det(m)!=0))),l->l!=\{\});}\\

This line \texttt{o6} is the set of all simple morphisms amongst lines \texttt{E\_0} to \texttt{E\_3}.\\

\texttt{i7: Ximg=flatten apply(simples,i->apply(i,j->apply(X,x->j*x)));}\\

This line \texttt{o7} is the set of all images of \texttt{X} under the simple morphisms.\\

\texttt{i8: all(Ximg,i->all(i,j->any(X,x->rank(x|j)==1)))}\\

We get \texttt{o8=true}, indicating that the images of \texttt{X} under all simple morphisms of \texttt{Xsimples} are all projectively equivalent to \texttt{X}. The next three lines verify the same for \texttt{Y}.\\

\texttt{i9: Ysimples= apply(\{0,1,2,3\},i->select(apply(\{0,1,2,3\},\\j->proj(E\_i,E\_j,E\_4)),m->det(m)!=0));}\\

\texttt{i10: Yimg=flatten apply(Ysimples,i->apply(i,j->apply(X,x->j*x)));}\\

\texttt{i11: all(Yimg,i->all(i,j->any(Y,y->rank(y|j)==1)))}\\

%
%
%
%
%
%
%
%
%

Note that since $\pi(L'',L',L)=\pi(L,L',L'')\inv$ for all distinct $L,L',L''\in\Lcal^*$, it suffices to show that $\pi(L,L',L_8)(X)=Y$ for all $L,L'\in\{L_0,L_2,L_4,L_6\}$, as \texttt{o11} confirms. Since \texttt{o8} and \texttt{o11} output \texttt{true}, we may conclude that the images of the Penrose points are Penrose points. Thus, we have $$\pi(L\cap\Pcal)=L''\cap \Pcal.$$

    
    Therefore by Corollary 1 we know that for all $L\in\Ecal$, we have $$\Aut_{(\Ecal,\Pi^*)}(L)\leq\Aut_{\Pcal\cap L}(L)\cong A_4.$$

    Now we must show that $\Aut_{\Pcal\cap L}(L)\leq\Aut_{(\Ecal,\Pi^*)}(L)$ for $L\in\Ecal$. The only nontrivial subgroups of $\PGL(2,\C)$ isomorphic to a proper subgroup of $A_4$ are $C_2$ and $C_3$. Thus, it suffices to find both an element of order 2 and an element of order 3 in $\Aut_{(\Ecal,\Pi^*)}(L)$. Observe $$\pi(L_4,L_6,L_0)\circ\pi(L_0,L_2,L_4)=\begin{pmatrix}
        t&1\\
        0&t^2
    \end{pmatrix},$$ which has order 3, so $C_2\not\cong\Aut_{(\Ecal,\Pi^*)}(L)$, and $$\pi(L_4,L_6,L_0)\circ\pi(L_0,L_2,L_4)\circ\pi(L_4,L_8,L_0)\circ\pi(L_0,L_2,L_4)=\begin{pmatrix}
        -1&t\\
        t&1
    \end{pmatrix},$$ which has order 2, so $C_3\not\cong\Aut_{(\Ecal,\Pi^*)}(L)$. Thus $$\Aut_{(\Ecal,\Pi^*)}(L)=\Aut_{L\cap\Pcal}(L)\cong A_4,$$ and in general, the group of the groupoid is $$\Aut_{(\Ecal,\Pi^*)}\cong A_4.$$
\end{proof}
\subsection{The Groupoid of the Penrose Spread}
Now we will move on to the configuration $\Lcal$ of all ten lines $L_0,\dots,L_9$ from the beginning of the section, and identify the automorphism group of the groupoid $(\Lcal,\Pi)$.
\begin{prop}
 Let $\Lcal$ be the set of all ten lines $L_0,\dots,L_9$, and consider the groupoid $(\Lcal,\Pi)$. Then $$\Aut_{(\Lcal,\Pi)}\cong S_4.$$
\end{prop}
Before we prove the proposition, we must first explore forty new points which we will deem ``quasi-Penrose" points, denoted by $\Qcal$.

For $i\neq 8,9$, the four quasi-Penrose points on $L_i$ are given by the $\P^1$-parameters $$\begin{pmatrix}
    1&1
\end{pmatrix},\,\begin{pmatrix}
    -1&1
\end{pmatrix},\,\begin{pmatrix}
    2t+1&1
\end{pmatrix},\,\begin{pmatrix}
    -1&2t+1
\end{pmatrix}.$$ For $L_8$ and $L_9$, the quasi-Penrose points are given by the $\P^1$-parameters $$\begin{pmatrix}
    t&1
\end{pmatrix},\,\begin{pmatrix}
    -t&1
\end{pmatrix},\,\begin{pmatrix}
    1&t-1
\end{pmatrix},\,\begin{pmatrix}
    3&t-1
\end{pmatrix}.$$

\begin{lemma}
    Let $$\wX=\left\{\begin{pmatrix}
        1\\0
    \end{pmatrix},\begin{pmatrix}
        0\\1
    \end{pmatrix},\begin{pmatrix}
        t\\1
    \end{pmatrix},\begin{pmatrix}
        -t^2\\1
    \end{pmatrix},\begin{pmatrix}
        1\\1
    \end{pmatrix},\begin{pmatrix}
        -1\\1
    \end{pmatrix},\begin{pmatrix}
        2t+1\\1
    \end{pmatrix},\begin{pmatrix}
        -1\\2t+1
    \end{pmatrix}\right\}\subseq\P^1_\C.$$ Then $\Aut_{\wX}(\P^1)\cong S_4$. In other words, the eight points of $\wX$ have octahedral (or cubic) symmetry.
\end{lemma}
\begin{proof}
    First note that $\Aut_{\wX}(\P^1_\C)\hookrightarrow S_8$, so $\Aut_{\wX}(\P^1_\C)$ is isomorphic to a subgroup of $S_8$. In particular, $\Aut_{\wX}(\P^1)$ must be finite.

    A straightforward analysis reveals that $\Aut_X(\P^1)\leq\Aut_{\wX}(\P^1)$; that is, every automorphism of $\P^1$ that permutes the points of $X$ also permutes the points of $\wX$. Thus $\Aut_{\wX}(\P^1)$ contains a subgroup isomorphic to $A_4$. Together with knowing that $\Aut_{\wX}(\P^1)$ is a finite subgroup of $\PGL(2,\C)$, we know that the only possibilities are $\Aut_{\wX}(\P^1)\cong A_4$ or $\Aut_{\wX}(\P^1)\cong S_4$.

    We can observe that the matrix $f=\begin{pmatrix}
        -1&2t+1\\
        1&1\\
    \end{pmatrix}\in\Aut_{\wX}(\P^1)\setminus\Aut_X(\P^1)$. Indeed, we have the table
\begin{center}
\begin{tabular}{c|c|c|c|c|c|c|c|c}
$x$
&
     $\begin{pmatrix}
        1\\0
    \end{pmatrix}$&$\begin{pmatrix}
        0\\1
    \end{pmatrix}$&$\begin{pmatrix}
        t\\1
    \end{pmatrix}$&$\begin{pmatrix}
        -t^2\\1
    \end{pmatrix}$&$\begin{pmatrix}
        1\\1
    \end{pmatrix}$&$\begin{pmatrix}
        -1\\1
    \end{pmatrix}$&$\begin{pmatrix}
        2t+1\\1
    \end{pmatrix}$&$\begin{pmatrix}
        -1\\2t+1
    \end{pmatrix}$ \\
    \hline
     $f(x)$
&
     $\begin{pmatrix}
        -1\\1
    \end{pmatrix}$&$\begin{pmatrix}
        2t+1\\1
    \end{pmatrix}$&$\begin{pmatrix}
        1\\1
    \end{pmatrix}$&$\begin{pmatrix}
        -1\\2t+1
    \end{pmatrix}$&$\begin{pmatrix}
        t\\1
    \end{pmatrix}$&$\begin{pmatrix}
        1\\0
    \end{pmatrix}$&$\begin{pmatrix}
        0\\1
    \end{pmatrix}$&$\begin{pmatrix}
        -t^2\\1
    \end{pmatrix}$ \\
\end{tabular} 
.
\end{center}
Thus $\Aut_X(\P^1)$ is a proper subgroup of $\Aut_{\wX}(\P^1)$ and therefore $\Aut_{\wX}(\P^1)\cong S_4$.
\end{proof}
In fact, we can write the entire set of $\Aut_{\wX}(\P^1)$ as \begin{align*}
\left\{
\begin{pmatrix}
    	1&0\\0&1
	\end{pmatrix},\begin{pmatrix}
    	0&-1\\1&0
	\end{pmatrix},\begin{pmatrix}
    	t&1\\1&-t
	\end{pmatrix},\begin{pmatrix}
    	-1&t\\t&1
	\end{pmatrix},\begin{pmatrix}
    	0&-1\\t&-t^2
	\end{pmatrix},\begin{pmatrix}
    	t&-t^2\\1&0
	\end{pmatrix},\right. \\\begin{pmatrix}
    	0&-t\\1&t^2
	\end{pmatrix},\begin{pmatrix}
    	t^2&t\\-1&0
	\end{pmatrix},\begin{pmatrix}
    	t&0\\1&t^2
	\end{pmatrix},\begin{pmatrix}
    	t^2&0\\-1&t
	\end{pmatrix},\begin{pmatrix}
    	t^2&t\\0&1
	\end{pmatrix},\begin{pmatrix}
    	t&-t^2\\0&1
	\end{pmatrix},\\
	\left(\begin{array}{cc}
-1&2t+1\\
1&1
\end{array}\right),\left(\begin{array}{cc}
-1&-1\\
-1&2t+1
\end{array}\right),\left(\begin{array}{cc}
2t+1&1\\
1&-1
\end{array}\right),\\\left(\begin{array}{cc}
-1&1\\
2t+1&1
\end{array}\right),\left(\begin{array}{cc}
-1&-1\\
1&-1
\end{array}\right),\left(\begin{array}{cc}
1&-1\\
-1&2t+1
\end{array}\right),\\\left(\begin{array}{cc}
-1&-1\\
2t+1&1
\end{array}\right),\left(\begin{array}{cc}
2t+1&1\\
1&-2t-1
\end{array}\right),\left(\begin{array}{cc}
-1&2t+1\\
2t+1&1
\end{array}\right),\\ \left.\left(\begin{array}{cc}
-1&2t+1\\
-1&1
\end{array}\right),\left(\begin{array}{cc}
2t+1&1\\
1&1
\end{array}\right),\left(\begin{array}{cc}
1&-1\\
1&1
\end{array}\right)\right\}
\end{align*}
\begin{lemma}
    Let $$\wY=\left\{\begin{pmatrix}
        1\\0
    \end{pmatrix},\begin{pmatrix}
        0\\1
    \end{pmatrix},\begin{pmatrix}
        -1\\1
    \end{pmatrix},\begin{pmatrix}
        t^2\\1
    \end{pmatrix},\begin{pmatrix}
        t\\1
    \end{pmatrix},\begin{pmatrix}
        -t\\1
    \end{pmatrix},\begin{pmatrix}
        1\\t-1
    \end{pmatrix},\begin{pmatrix}
        3\\t-1
    \end{pmatrix}\right\}\subseq\P^1_\C.$$ Then $\Aut_{\wY}(\P^1)\cong S_4$.
\end{lemma}
\begin{proof}
    The proof is identical to that of Lemma 3, except the automorphism group $\Aut_{\wY}(\P^1)$ is equal to the set
    \begin{align*}
    \left\{\begin{pmatrix}
        1&0\\0&1
    \end{pmatrix},\begin{pmatrix}
        0&-1\\t&0
    \end{pmatrix},\begin{pmatrix}
        -1&t^2\\1&1
    \end{pmatrix},\begin{pmatrix}
        1&1\\t&-1
    \end{pmatrix},\begin{pmatrix}
        0&-1\\1&1
    \end{pmatrix},\begin{pmatrix}
        1&1\\-1&0
    \end{pmatrix}, \right.\\\begin{pmatrix}
        0&t\\-1&t^2
    \end{pmatrix},\begin{pmatrix}
        t^2&-t\\1&0
    \end{pmatrix},\begin{pmatrix}
        t&0\\-t&1
    \end{pmatrix},\begin{pmatrix}
        1&0\\t&t
    \end{pmatrix},\begin{pmatrix}
        t&-1\\0&1
    \end{pmatrix},\begin{pmatrix}
        1&1\\0&t
    \end{pmatrix},\\
        \left(\begin{array}{cc} 1&t\\ t-1&-1 \end{array}\right),\left(\begin{array}{cc} t-1&-1\\ -t&t+1 \end{array}\right),\left(\begin{array}{cc} t+1&1\\ t&t-1 \end{array}\right),\\\left(\begin{array}{cc} t&t-1\\ 1&-t \end{array}\right),\left(\begin{array}{cc} -t+1&1\\ t&t-1 \end{array}\right),\left(\begin{array}{cc} t&t-1\\ -1&-t \end{array}\right),\\\left(\begin{array}{cc} -2t-1&-t\\ t+1&1 \end{array}\right),\left(\begin{array}{cc} t&t+1\\ 1&t \end{array}\right),\left(\begin{array}{cc} t&-t-1\\ -1&t \end{array}\right),\\\left.\left(\begin{array}{cc} 1&t\\ -t-1&-2t-1 \end{array}\right),\left(\begin{array}{cc} -t-1&-1\\ t+2&t+1 \end{array}\right),\left(\begin{array}{cc} t&t-1\\ -2t-1&-t \end{array}\right)\right\}.\\
    \end{align*}
\end{proof}
\begin{cor}
    Let $L\in\Lcal$ and consider the Penrose and quasi-Penrose points $\Pcal\cup \Qcal$ on $L$. Then $$\Aut_{L\cap(\Pcal\cup\Qcal)}(L)\cong S_4.$$
\end{cor}
\begin{proof}
    Follows immediately from Lemmas 3 and 4.
\end{proof}
\begin{thm}\label{PS4}
    Let $\Pi$ be the composition-closure of the set of morphisms of the form $\pi(L,L',L'')$, with $L,L',L''\in\Lcal$ distinct. Then the groupoid $(\Lcal,\Pi)$ has the automorphism group $$\Aut_{(\Lcal,\Pi)}\cong S_4.$$
\end{thm}
\begin{proof}
    We will show that for any $P\in\Pcal\cup \Qcal$ and $\pi\in\Pi$, we have $$\pi(P)\in\Pcal\cup \Qcal.$$ We will begin by proving the above claim in the case $\pi$ is a simple morphism of $\Pi$, using the following Macaualy2 code \cite{M2}.\\
    
\texttt{i1: K=QQ[t]/ideal(t\^{}2+t+1)}\\

    \texttt{i2: L=\{\{\{1,0,t+1,t+1\},\{1,0,-t,-t\}\},\{\{0,1,t+1,-t-1\},\{0,1,-t,t\}\},\\\{\{1,0,t+1,-t\},\{0,1,-t,t\^{}2\}\},\{\{0,1,-t\^{}2,t\},\{1,0,-t,-t\^{}2\}\},\\\{\{0,1,-t,1\},\{1,0,-1,-t\}\},\{\{1,0,-t,-1\},\{0,1,-1,t\}\},\\\{\{1,0,-1,-t\^{}2\},\{0,1,-t\^{}2,1\}\},\{\{0,1,-1,t\^{}2\},\{1,0,-t\^{}2,-1\}\},\\\{\{1,1,t,0\},\{1,1,t\^{}2,0\}\},\{\{1,-1,0,t\},\{1,-1,0,t\^{}2\}\}\};}\\

\texttt{i3: proj=(i,j,k)->matrix\{\{-det(matrix\{i\_0,j\_0,j\_1,k\_1\}),\\-det(matrix\{i\_1,j\_0,j\_1,k\_1\})\},\{det(matrix\{i\_0,j\_0,j\_1,k\_0\}),\\det(matrix\{i\_1,j\_0,j\_1,k\_0\})\}\};}\\

\texttt{i4: X=\{matrix\{\{1\},\{0\}\},matrix\{\{0\},\{1\}\},matrix\{\{t\},\{1\}\},\\matrix\{\{-t\^{}2\},\{1\}\},matrix\{\{1\},\{1\}\},matrix\{\{-1\},\{1\}\},\\matrix\{\{2*t+1\},\{1\}\},matrix\{\{-1\},\{2*t+1\}\}\};}\\

\texttt{i5: Y=\{matrix\{\{1\},\{0\}\},matrix\{\{0\},\{1\}\},matrix\{\{-1\},\{1\}\},\\matrix\{\{t\^{}2\},\{1\}\},matrix\{\{t\},\{1\}\},matrix\{\{-t\},\{1\}\},\\matrix\{\{1\},\{t-1\}\},matrix\{\{3\},\{t-1\}\}\};}\\

\texttt{i6: simples=(a,b)->apply((select(apply(L,l->proj(L\_a,l,L\_b)),\\i->det(i)!=0)),i->apply(X,x->i*x));}\\

This command in \texttt{i6} receives an input of two indices $a$ and $b$ and computes $\pi(L_a,\ell,L_b)(\wX)$ as lists for all $\ell\in \Lcal$, where $\det(\pi(L_a,\ell,L_b))\neq0$.\\

\texttt{i7: all(apply(subsets(toList(\{0..7\}\_0),2),r->all(simples(r\_0,r\_1),\\s->all(s,i->any(X,x->det(x|i)==0)))),w->w==true)}\\

This command \texttt{i7} takes each of the sets \texttt{simples(r\_0,r\_1)} for $r_0,r_1\in\{0,\dots,7\}$ and verifies whether each set within \texttt{simples(r\_0,r\_1)} is projectively equivalent to $\wX$ by going through each member in each set contained in \texttt{simple(r\_0,r\_1)} and determining whether there exists an element of $\wX$ equivalent to that member by joining the two $2\times 1$ matrices together (one from \texttt{simples(r\_0,r\_1)\_i} and one from $\wX$) into one $2\times 2$ matrix and determining whether the determinant of that $2\times 2$ matrix is 0. If every element of every set inside \texttt{simple(r\_0,r\_1)} is equivalent to some member of $\wX$, then the image of $\wX$ under the invertible map $\pi(L_{r_1},\ell,L_{r_2})$ is $\wX$ for every $\ell\in \Lcal$.\\ 

\texttt{i8: all(apply(toList(\{0..7\}\_0),r->all(simples(r,8),\\s->all(s,i->any(Y,y->det(y|i)==0)))),w->w==true)}\\

\texttt{i9:all(apply(toList(\{0..7\}\_0),r->all(simples(r,9),\\s->all(s,i->any(Y,y->det(y|i)==0)))),w->w==true)}\\

The above commands \texttt{i8} and \texttt{i9} verify each of the invertible simple morphisms $\pi(L_r,\ell, L_8)$ and $\pi(L_r,\ell,L_9)$ map $\wX$ to $\wY$, respectively, for $L_r,\ell\in\Lcal$.\\

\texttt{i10: apply((select(apply(L,l->proj(L\_8,l,L\_9)),\\i->det(i)!=0)),i->apply(Y,j->i*j));}\\

This command \texttt{i10} creates a list of images of $\wY$ under each of the invertible simple morphisms $\pi(L_8,\ell,L_9)$ for $\ell\in\Lcal$.\\

\texttt{i11: all(o10,s->all(s,i->any(Y,y->det(y|i)==0)))}\\

The above command \texttt{i11} verifies each of the invertible simple morphisms $\pi(L_8,\ell, L_9)$ map $\wY$ to $\wY$, for $\ell\in\Lcal$.\\

The values of \texttt{o7}, \texttt{o8}, \texttt{o9}, and \texttt{o11} are all \texttt{true}, meaning that for any simple morphism $\pi(L,L',L'')$, we have $$\pi(L,L',L'')(L\cap(\Pcal\cup\Qcal))=L''\cap(\Pcal\cup \Qcal).$$ Since $\Pi$ is the composition-closure of simple morphisms, it holds that for all $(\pi:L\to L'')\in\Pi$, we have $$\pi(L\cap(\Pcal\cup\Qcal))=L''\cap(\Pcal\cup\Qcal).$$ Therefore by Corollary 2 we know that for all $L\in\Lcal$, we have $$\Aut_{(\Lcal,\Pi)}(L)\leq\Aut_{L\cap(\Pcal\cup\Qcal)}(L)\cong S_4.$$ 

Now we must show that $\Aut_{L\cap (\Pcal\cup\Qcal)}(L)\leq \Aut_{(\Lcal,\Pi)}(L)$ for $L\in\Lcal$. Note that since $(\Ecal,\Pi^*)$ is a subgroupoid of $(\Lcal,\Pi)$, we have for all $L\in\Ecal\subseq\Lcal$, $$A_4\cong\Aut_{L\cap\Pcal}(L)=\Aut_{(\Ecal,\Pi^*)}(L)\leq\Aut_{(\Lcal,\Pi)}(L)\leq\Aut_{L\cap(\Pcal\cup\Qcal)}(L)\cong S_4.$$ So to prove $\Aut_{L\cap\Pcal}(L)$ is a proper subgroup of $\Aut_{(\Lcal,\Pi)}(L)$, it is sufficient to find in $\Aut_{(\Lcal,\Pi)}(L)$ an automorphism that permutes $L\cap(\Pcal\cup\Qcal)$ and does not permute $L\cap \Pcal$. An example of such an automorphism is $$\pi(L_2,L_4,L_0)\circ\pi(L_0,L_1,L_2)=\begin{pmatrix}
    -1&1\\
    1&-2t-1
\end{pmatrix}\in\Aut_{L_0\cap(\Pcal\cup\Qcal)}(L_0)\setminus\Aut_{L_0\cap\Pcal}(L_0).$$ Thus for $L\in\Lcal$ we have $\Aut_{L\cap\Pcal}(L)$ is a proper subgroup of $\Aut_{(\Lcal,\Pi)}(L)$, and so $$\Aut_{(\Lcal,\Pi)}(L)=\Aut_{L\cap(\Pcal\cup\Qcal)}(L)\cong S_4,$$ and so in general $$\Aut_{(\Lcal,\Pi)}\cong S_4.$$
\end{proof}
\subsection{The Double Penrose Configuration}\
This section shall study the union of the Penrose points $\Pcal$ and the quasi-Penrose points $\Qcal$ that have been illuminated by the groupoid $(\Lcal,\Pi)$ from the previous section. We will call this union the double-Penrose configuration, denoted $\Dcal$. First note that $\Qcal$ is in fact isomorphic to $\Pcal$ as sets of points $\P^3_\C$; indeed, a computation in Macaulay2 \cite{M2} verifies that $$\Qcal^{\mathsf{T}}=\underbrace{\begin{pmatrix}
0&0&1&1\\
0&0&1&-1\\
-1&-1&0&0\\
-1&1&0&0\\
\end{pmatrix}}_{A}\Pcal^{\mathsf{T}},$$ where $A\in\PGL(4,\C)=\Aut(\P^3_\C)$.

Therefore we know that $\Qcal$ must be a $(5,8)$-geproci non-half-grid. Experiments in Macaulay2 projecting $\Dcal$ from a random point of $P\in\P^3_\Q$ have found octic cones with a singularity of order eight at $P$. This indicates that $\Dcal$ is almost certainly $(8,10)$-geproci, but a formal proof eludes us, as projecting from a general point has been too computationally expensive. If $\Dcal$ is $(8,10)$-geproci, it would be a half-grid because the degree-10 curve containing the projection of $\Dcal$ can be taken to be the projection of the ten lines $L_0,\dots L_9$.
We have the following ideals at our disposal. $$I_{\Pcal}=\left(x\,y\,z\,w,\,x^{3}w+y^{3}w+z^{3}w,\,x^{3}z-y^{3}z+
       z\,w^{3},\,x\,y^{3}-x\,z^{3}+x\,w^{3},\,x^{3}y-y\,z^{3}-y\,w^{3}\right),$$
       and
       \begin{align*}I_{\Qcal}=(x^{2}z^{2}-y^{2}z^{2}-x^{2}w^{2}+y^{2}w^{2},\\\,x^{4
       }+2\,x^{3}y-2\,x\,y^{3}-y^{4}-2\,x\,z^{3}+2\,y\,z^{3}-6\,x\,z\,w^{2}+6\,y
       \,z\,w^{2},\\ x^{4}-2\,x^{3}y+2\,x\,y^{3}-y^{4}-6\,x\,z^{2}w-6\,y\,z^{2}w-
       2\,x\,w^{3}-2\,y\,w^{3},\\\,6\,x^{2}y\,z+2\,y^{3}z+z^{4}+6\,x^{2}y\,w+2\,y
       ^{3}w-2\,z^{3}w+2\,z\,w^{3}-w^{4},\\\,2\,x^{3}z+6\,x\,y^{2}z+z^{4}-2\,x^{3}
       w-6\,x\,y^{2}w+2\,z^{3}w-2\,z\,w^{3}-w^{4}),\end{align*}
and thus we can construct $I_\Dcal=I_\Pcal\cap I_\Qcal$. We can test whether $\Dcal$ is geproci by moving up to the ring $R=\C[x,y,z,w][a,b,c,d]$ and defining $P=(dx-aw,dy-bw,dz-cw)$, and testing whether the intersection $(I_{\Dcal}*R)\cap P^8$ contains any octics. So far, this has proven to be too computationally expensive, but perhaps further study will reveal a less computationally-heavy proof.
\section{The Klein Configuration}
Let $i\in \C$ be the imaginary unit $i^2=-1$. Then the Klein configuration is defined as the following set $\Kcal$ of sixty points of $\P^3_\C$.
\begin{center}
\begin{tabular}{|c|c|c|c|}
\hline
 $(0,0,1,1)$ &  $(0,0,1,i)$ &  $(0,0,1,-1)$ &  $(0,0,1,-i)$\\
\hline
 $(0,1,0,1)$ &  $(0,1,0,i)$ &  $(0,1,0,-1)$ & $(0,1,0,-i)$ \\
\hline
 $(0,1,1,0)$ &  $(0,1,i,0)$ & $(0,1,-1,0)$ &  $(0,1,-i,0)$ \\
\hline
 $(1,0,0,1)$ &  $(1,0,0,i)$ &  $(1,0,0,-1)$ &  $(1,0,0,-i)$ \\
\hline
 $(1,0,1,0)$ &  $(1,0,i,0)$ &  $(1,0,-1,0)$ &  $(1,0,-i,0)$ \\
\hline
 $(1,1,0,0)$ &  $(1,i,0,0)$ &  $(1,-1,0,0)$ &  $(1,-i,0,0)$ \\
\hline
$(1,0,0,0)$ &  $(0,1,0,0)$ & $(0,0,1,0)$  & $(0,0,0,1)$ \\
\hline
 $(1,1,1,1)$ &  $(1,1,1,-1)$ &  $(1,1,-1,1)$ & $(1,1,-1,-1)$ \\
\hline
 $(1,-1,1,1)$ &  $(1,-1,1,-1)$ & $(1,-1,-1,1)$ & $(1,-1,-1,-1)$ \\
\hline
$(1,1,i,i)$ & $(1,1,i,-i)$ &  $(1,1,-i,i)$ & $(1,1,-i,-i)$ \\
\hline
 $(1,-1,i,i)$ &  $(1,-1,i,-i)$ &  $(1,-1,-i,1)$ &  $(1,-1,-i,-i)$ \\
\hline
$(1,i,1,i)$ &  $(1,i,1,-i)$ &  $(1,-i,1,i)$ &  $(1,-i,1,-i)$ \\
\hline
 $(1,i,-1,i)$ & $(1,i,-1,-i)$ &  $(1,-i,-1,i)$ & $(1,-i,-1,-i)$ \\
\hline
 $(1,i,i,1)$ & $(1,i,-i,1)$ &  $(1,-i,i,1)$ &  $(1,-i,-i,1)$ \\
\hline
 $(1,i,i,-1)$ &  $(1,i,-i,-1)$ & $(1,-i,i,-1)$ &  $(1,-i,-i,-1)$\\
 \hline
 \end{tabular}
 \end{center}

The Klein configuration is known to be a $(6,10)$-geproci non-half-grid \cite{PSS}. Furthermore, $\Kcal$ is a $(60_3,30_6)$-configuration of points and lines. We may define a set $\Lcal$ of ten mutually-skew lines $L_0,\dots,L_9$ defined as the rowspans of the respective matrices $\bL_i$, given as follows.
\begin{align*}
    \bL_0&=\left(\begin{array}{cccc}
       0&0&1&0\\
       0&0&0&1
       \end{array}\right)&\bL_1&=\left(\begin{array}{cccc}
       0&1&0&-i\\
       1&0&i&0
       \end{array}\right)\\\bL_2&=\left(\begin{array}{cccc}
       0&1&1&0\\
       1&0&0&-1
       \end{array}\right)&\bL_3&=\left(\begin{array}{cccc}
       0&1&i&0\\
       1&0&0&i
       \end{array}\right)\\\bL_4&=\left(\begin{array}{cccc}
       0&1&-1&0\\
       1&0&0&1
       \end{array}\right)&\bL_5&=\left(\begin{array}{cccc}
       0&1&-i&0\\
       1&0&0&-i
       \end{array}\right)\\\bL_6&=\left(\begin{array}{cccc}
       0&1&0&1\\
       1&0&1&0
       \end{array}\right)&\bL_7&=\left(\begin{array}{cccc}
       1&0&0&0\\
       0&1&0&0
       \end{array}\right)\\\bL_8&=\left(\begin{array}{cccc}
       0&1&0&i\\
       1&0&-i&0
       \end{array}\right)&\bL_9&=\left(\begin{array}{cccc}
       0&1&0&-1\\
       1&0&-1&0
       \end{array}\right)
       \end{align*}


We claim that the group of the groupoid is isomorphic to $S_4$, and the six Klein points on each line have octahedral symmetry.

One can verify that each of the 60 points of $\Kcal$ is equal to a linear combination of the rows of one of the matrices $\bL_i$ with the coefficients given by the vectors of $$E=\left\{\begin{pmatrix}
    1\\0
\end{pmatrix}\begin{pmatrix}
    0\\1
\end{pmatrix}\begin{pmatrix}
    1\\1
\end{pmatrix}\begin{pmatrix}
    -1\\1
\end{pmatrix}\begin{pmatrix}
    i\\1
\end{pmatrix}\begin{pmatrix}
    -i\\1
\end{pmatrix}\right\}.$$ That is, $$\Kcal=\left\{\vec{e}^{\mathsf{T}}\cdot\bL_i:\vec{e}\in E,0\leq i\leq 9\right\}.$$

\begin{lemma}\label{Kleinocta}
    Let $G=\Aut_E(\P^2_\C)$ be the subgroup of $\Aut(\P^1_\C)$ that setwise fixes $E$. Then $G\cong S_4$.
\end{lemma}
\begin{proof}
The set of points $E\subseq\widehat{\C}$ forms the vertices of an octahedron in the Riemann sphere. The set of rotations of this set is the octahedral group, isomorphic to $S_4$.

    In more detail, first note that since $G$ acts on a set of six points, and the images of three distinct points uniquely determine the identity of an automorphism of $\P^1$, that $G\hookrightarrow S_6$. Specifically, we know $G$ must be finite.

    Note that the two automorphisms $\begin{pmatrix}
        i&0\\
        0&1\\
    \end{pmatrix},\begin{pmatrix}
        1&-1\\1&1
    \end{pmatrix}\in G$ can be embedded into the group $\text{SO}(3)$ of rotations of $\R^3$ given by given by $\begin{pmatrix}i&0\\0&1\end{pmatrix}\mapsto\begin{pmatrix}0&1&0\\-1&0&0\\0&0&1\\\end{pmatrix}$ and $\begin{pmatrix}1&-1\\1&1\end{pmatrix}\mapsto\begin{pmatrix}0&0&1\\0&1&0\\-1&0&0\\\end{pmatrix}$, by embedding the Riemann sphere into $\R^3$. According to GAP, the group generated by these two $3\times 3$ matrices is isomorphic to the symmetric group $S_4$ \cite{GAP}. Up to isomorphism, only finite subgroup of $\PGL(2,\C)\cong\Aut(\P^1_\C)$ containing a subgroup isomorphic to $S_4$ is $S_4$ itself. Thus $G\cong S_4$.
    

    
%
\end{proof}
\begin{cor}\label{Kleinary}
    Let $L\in\Lcal$. Then $$\Aut_{L\cap\Kcal}(L)\cong S_4.$$
\end{cor}
\begin{proof}
    Follows directly from Lemma \ref{Kleinocta}.
\end{proof}
\begin{thm}\label{KS4}
    Let $\Pi$ be the composition-closure of the set of morphisms of the form $\pi(L,L',L'')$, with $L,L',L''\in\Lcal$. Then the groupoid $(\Lcal,\Pi)$ has the automorphism group $$\Aut_{(\Lcal,\Pi)}\cong S_4.$$
\end{thm}
\begin{proof}
    We will show that for any $P\in \Kcal$ and $\pi\in \Pi$, that $\pi(P)\in \Kcal$. We will begin by proving the above claim in the case that $\pi$ is a simple morphism of $\Pi$, using the following Macaulay2 code \cite{M2}. \\

    \texttt{i1: F=QQ[i]/ideal(i\^{}2+1);}\\
    
    \texttt{i2: L=\{\{\{0, 0, 1, 0\}, \{0, 0, 0, 1\}\}, \{\{0, 1, 0, -i\}, \{1, 0, i, 0\}\},\\ \{\{0, 1, 1,0\}, \{1, 0, 0, -1\}\}, \{\{0, 1, i, 0\}, \{1, 0, 0, i\}\},\\ \{\{0, 1, -1, 0\}, \{1, 0,
       0, 1\}\}, \{\{0, 1, -i, 0\}, \{1, 0, 0, -i\}\},\\ \{\{0, 1, 0, 1\}, \{1, 0, 1, 0\}\},
       \{\{1, 0, 0, 0\}, \{0, 1, 0, 0\}\},\\ \{\{0, 1, 0, i\}, \{1, 0, -i, 0\}\}, \{\{0, 1, 0,
       -1\}, \{1, 0, -1, 0\}\}\};}\\

       \texttt{i3: proj=(i,j,k)->matrix\{\{-det(matrix\{i\_0,j\_0,j\_1,k\_1\}),\\-det(matrix\{i\_1,j\_0,j\_1,k\_1\})\},\{det(matrix\{i\_0,j\_0,j\_1,k\_0\}),\\det(matrix\{i\_1,j\_0,j\_1,k\_0\})\}\};}\\

       \texttt{i4: E=\{matrix\{\{1\},\{0\}\}, matrix\{\{0\},\{1\}\}, matrix\{\{1\},\{1\}\}, matrix\{\{-1\},\{1\}\},\\ matrix\{\{i\},\{1\}\}, matrix\{\{-i\},\{1\}\}\};}\\

       \texttt{i5: allSimples=unique flatten flatten apply(L,t->apply(select(L,u->u!=t),u->\\apply(select(L,v->v!=u and v!=t),v->proj(t,u,v))));}\\

       This command in \texttt{i5} computes all simple morphisms of $\Pi$ and puts them together in one list. (Note this list will still contain repeated matrices up to scalar-invariance in $\PGL(2,\C)$; this does not interfere with the argument, so we did not remove the repeated morphisms.)\\

       \texttt{i6: allOutputs=unique flatten apply(E,e->apply(allSimples,m->m*e));}\\

       This command in \texttt{i6} computes $\{\pi(E):\pi\in\Pi\text{ is simple}\}$.\\

       \texttt{i7: all(allOutputs,j->any(E,e->rank(j|e)==1))}\\

       This command in \texttt{i7} verifies that for every simple morphism $\pi$ and $\vec{e}\in E$, $\pi(\vec{e})\in E$. 
       
       Thus, with $\pi$ viewed as an automorphism of $\P^1_\C$, we have $\pi(E)=E$ for every simple morphism $\pi$. Therefore $$\pi(L,L',L'')(\Kcal\cap L)=\Kcal\cap L''.$$ Since $\Pi$ is the composition-closure of simple morphims, it holds that for all $(\pi:L\to L'')\in\Pi$, we have $$\pi(\Kcal\cap L)=\Kcal\cap L''.$$ Therefore by Corollary \ref{Kleinary}, we know that for all $L\in \Lcal$, we have $$\Aut_{(\Lcal,\Pi)}(L)\leq\Aut_{\Kcal\cap L}(L)\cong S_4.$$ Now we must show that $\Aut_{\Kcal\cap L}(L)\leq \Aut_{(\Lcal,\Pi)}(L)$. 
       
       The only nontrivial subgroups of $\PGL(2,\C)$ isomorphic to a subgroup of $S_4$ are isomorphic to $C_2$, $C_3$, $C_4$, $D_4$, $A_4$, or $S_4$. The examples $$\pi(L_2,L_3,L_0)\circ\pi(L_0,L_1,L_2)=\begin{pmatrix}
           1&-i\\
           1&i
       \end{pmatrix}$$ and $$\pi(L_2,L_7,L_0)\circ\pi(L_0,L_1,L_2)=\begin{pmatrix}
           1&-i\\
           -i&1
       \end{pmatrix}$$ show that $\Aut_{(\Lcal,\Pi)}(L_0)$ contains elements of orders $3$ and $4$, and so out of the six options, we must have $\Aut_{(\Lcal,\Pi)}(L_0)\cong S_4$ \cite{groupnames}. Thus $\Aut_{(\Lcal,\Pi)}\cong S_4$.

\end{proof}
\section{Future Work}
\subsection{Groupoids with Infinite and Icosahedral Automorphism Group}
Experiments in Macaulay2 have revealed not all configurations connected to the geproci property induce groupoids with finite automorphism group. Indeed, testing the lines spanning the points of the $F_4$ configuration, which is a $(4,6)$-geproci configuration \cite{politus1}, has revealed group elements with infinite order, so the group itself must be infinite.

Furthermore, we have looked at the 120-cell described in Remark 8.2 of \cite{politus1} to gather data on the groupoid it induces. The 120-cell is a $(10,12)$-geproci non-half-grid arranged in a $(120_4,80_6)$-configuration of points and lines. We have found a set $\Lcal$ of 20 mutually-disjoint lines that span all 120 points of the configuration; experiments in Macaulay2 have revealed group elements of infinite order, and so the automorphism group $\Aut_{(\Lcal,\Pi)}$ must be infinite.

 Similar experiments have revealed group elements of infinite order induced by the 27 lines of the cubic surface $V(x^3+y^3+z^3+w^3-(x+y+z+w)^3)$, and the 12 lines of the Sch\"afli Double Six. In the future, it would be interesting to classify these infinite automorphism groups, and to find interesting geometric arrangements thereby induced, such as orbits of finite size.
 
We are also interested in finding a configuration of lines $\Lcal\subseteq\P^3_{\C}$ whose group of the groupoid $\Aut_{(\Lcal,\Pi)}$ is isomorphic to the icosohedral group $A_5$. This is the largest finite subgroup of $\PGL(2,\C)$ \cite{Klein}, and is the only alternating group contained in $\PGL(2,\C)$ we have not yet found.
 
\subsection{Higher-Dimensional Configurations}
Let us use Construction 4.2 in \cite{Dolgachev} with $p=5$ to get a configuration of 25 lines in $\P^4_\C$. Let $t\in\C$ satisfy $t^4+t^3+t^2+t+1=0$.  Let $G\leq \Aut(\P^4_\C)$ be generated by the automorphims $$\sigma(p_0,p_1,p_2,p_3,p_4)=(p_4,p_0,p_1,p_2,p_3)$$ and $$\tau(p_0,p_1,p_2,p_3,p_4)=(p_0,tp_1,t^2p_2,t^3p_3,t^4p_4).$$ Then $G\cong C_5^2$, and the orbit of the line $$\mathbf{L}_{0,0}=\begin{pmatrix}0&1&0&0&-1\\
0&0&1&-1&0\\\end{pmatrix}$$ forms a configuration of 25 lines $\Lcal_{25}$. The lines form a $(25_6,30_5)$-configuration with a set of 30 three-dimensional hyperplanes \cite{Dolgachev}.  
The ideal $I$ of the 25 lines is minimally generated by \begin{align*}
    xyzuv,\\
    x^3zu+xyu^3+xz^3v+y^3uv+yzv^3,\\
    x^2y^2u+xz^2u^2+y^2z^2v+x^2zv^2+yu^2v^2,\\
    x^2yz^2+y^2zu^2+x^2u^2v+xy^2v^2+z^2uv^2,\\
    xy^3z+yz^3u+x^3yv+zu^3v+xuv^3,\\
    x^5+y^5+z^5+u^5+v^5.
\end{align*}

Experiments in Macaulay2 indicate that $\Lcal_{25}$ is almost certainly not geproci; the ideal of its projection from a random rational point into $\P^3_\C$ have been minimally generated by twenty septics and forty-two octics.

%
%

Despite the arrangement existing in $\P^4$, we can still define a groupoid induced by $\Lcal_{25}$ by using the 30 three-dimensional hyperplanes of the configuration. Let $U,V,W\in \Lcal_{25}$ and let $T\subseq \P^4$ be a three-dimensional hyperplane such that $U,V,W\subseq T$. Then we may define the projection isomorphism $\pi(U,V,W):U\to W$ by defining $\pi(U,V,W)(u)$ to be the unique point where the line $W$ meets the plane $\overline{uV}$. This intersection is nonempty because all three lines are contained in the three-dimensional hyperplane $T$.

This construction allows us to build a groupoid $(\Lcal_{25},\Pi)$ from this configuration of lines, where $\Pi$ is the composition-closure of maps of the form $\pi(U,V,W)$, such that there exists a three-dimensional hyperplane $T$ where $U,V,W\subseq T$.

Experiments in Macaulay2 have discovered group elements $g\in\Aut_{(\Lcal_{25},\Pi)}$ of infinite order. Interestingly, one such discovered group element is $\begin{pmatrix}1&1\\0&1\\\end{pmatrix}\in\PGL(2,\C)$, suggesting a connection to the Fibonacci numbers. In the future, it would be interesting to completely characterize $\Aut_{(\Lcal_{25},\Pi)}$ and to find new configurations induced by the orbits of points under the groupoid. In addition, we would also like to generalize the pattern exhibited by the minimal generators of $I$ to study yet higher-degree and -dimensional configurations.

\subsection{Configurations in Positive Characteristic}
Much of the groundwork for the study of arrangements of lines has been laid in the positive characteristic setting \cite{Kettinger, Ganger}. In fact, the connection between the groupoid induced by a spread of $\P^3_{\F_q}$ has been studied extensively \cite{Ganger}. In the future, we would like to study the groupoids induced by maximal partial spreads, which also have a well-documented connection to the geproci property \cite{Kettinger}. We would also like to find a generalizable method of finding maximal partial spreads, akin to the method used by Gorla, Managaniello, and Rosenthal to create full spreads \cite{Gorla}. Furthermore, we would like to examine the configurations studied in this paper in the positive characteristic setting; for example, it is known that $\P^3_{\F_7}$ contains a set of 40 points isomorphic to the Penrose configuration \cite{Kettinger}. Does the groupoid induced over $\F_7$ have the same automorphism group as that induced over $\C$?
\subsection{Lifts to the Grassmannian}
Pursuing an idea first proposed by Chris Peterson, we would like to generalize the notion of projection from a point (or line, plane, etc.) to a plane (or line, hyperplane, etc.) in projective space into a notion of projection from one Schubert cycle into another Schubert cycle in a Grassmannian, thereby generalizing the concept of geproci configurations.

For example, the arrangement $\Lcal_{25}$ as described in Section 6.2 can be expressed as an arrangement of 25 points in $\Gr(2,5)$. One may define a projection of $\lambda\in\Gr(2,5)$ from the general point $\gamma\in\Gr(2,5)$ into a Schubert $\Sigma_{2,2}$-cycle $\Pi$ by finding the point of intersection of $\Pi$ with the unique Schubert $\Sigma_{1,1}$-cycle $T$ containing both $\gamma$ and $\lambda$. This is well-defined because $\Sigma_{1,1}$ is the class of lines contained in a given three-dimensional hyperplane of $\P^4$; two general lines in $\P^4$ will be contained in a unique three-dimensional hyperplane, so $T$ will be unique. Furthermore, one may use Schubert calculus to verify that a general $\Sigma_{1,1}$ will intersect a $\Sigma_{2,2}$ at exactly a point, so the projection will be well-defined.

In an experiment in Macaulay, we have found that projecting the 25 points induced by $\Lcal_{25}$ from a random point $\rho\in\Gr(2,5)$ has not given a complete intersection, but rather a configuration whose ideal is generated by three sextics and two septics. This suggests with probability 1 that the 25 points are not geproci in this way; we would like to find arrangements that are.


\printbibliography

@book{G,
  title={Configurations of Points and Lines},
  author={Gr{\"u}nbaum, Branko},
  journal={Graduate Studies in Mathematics},
  volume={103},
  number={},
  pages={},
  year={2009},
  publisher={American Mathematical Society}
}

@book{Reye1876,
  title={Geometrie der Lage},
  author={Reye, Theodor},
  journal={},
  volume={2},
  number={1},
  pages={},
  year={1876},
  publisher={}
}

@book{Kantor1879,
  title={{\"U}ber eine Gattung von Configurationen in der Ebene und im Raume},
  author={Kantor, Seligmann},
  journal={Wien. Ber},
  volume={LXXXIV},
  number={},
  pages={915--932},
  year={1876},
  publisher={}
}

@online{ilya,
  title={Desargues Configurations and Groupoids},
  author={Zakharevich, Ilya},
  journal={},
  volume={},
  number={},
  pages={},
  year={2008},
  publisher={}
}

@article{Ganger,
  title={Spreads and Transversals and Their Connection to Geproci Sets},
  author={Ganger, Allison},
  journal={Dissertations, Theses, and Student Research Papers in Mathematics},
  volume={143},
  number={},
  pages={},
  year={2024},
  publisher={DigitalCommons@University of Nebraska - Lincoln}
}

@article{PSS,
author = {Piotr Pokora and Tomasz Szemberg and Justyna Szpond},
title = {{Unexpected Properties of the Klein Configuration of 60 Points in $\mathbb{P}^3$}},
volume = {74},
journal = {Michigan Mathematical Journal},
number = {3},
publisher = {University of Michigan, Department of Mathematics},
pages = {599 -- 615},
year = {2024},
doi = {10.1307/mmj/20216141},
URL = {https://doi.org/10.1307/mmj/20216141}
}

@misc{politusgroupoid,
      title={Combinatorics of skew lines in $\mathbb P^3$ with an application to algebraic geometry}, 
      author={Luca Chiantini and Łucja Farnik and Giuseppe Favacchio and Brian Harbourne and Juan Migliore and Tomasz Szemberg and Justyna Szpond},
      year={2025},
      eprint={2308.00761},
      archivePrefix={arXiv},
      primaryClass={math.AG},
      url={https://arxiv.org/abs/2308.00761}, 
}

@misc{politus1,
      title={Configurations of points in projective space and their projections}, 
      author={Luca Chiantini and Łucja Farnik and Giuseppe Favacchio and Brian Harbourne and Juan Migliore and Tomasz Szemberg and Justyna Szpond},
      year={2022},
      eprint={2209.04820},
      archivePrefix={arXiv},
      primaryClass={math.AG},
      url={https://arxiv.org/abs/2209.04820}, 
}

@misc{politus2,
      title={On the classification of certain geproci sets}, 
      author={Luca Chiantini and Lucja Farnik and Giuseppe Favacchio and Brian Harbourne and Juan Migliore and Tomasz Szemberg and Justyna Szpond},
      year={2024},
      eprint={2303.16263},
      archivePrefix={arXiv},
      primaryClass={math.AG},
      url={https://arxiv.org/abs/2303.16263}, 
}

@misc{politus3,
      title={Geproci sets on skew lines in $\mathbb P^3$ with two transversals}, 
      author={Luca Chiantini and Pietro De Poi and Lucja Farnik and Giuseppe Favacchio and Brian Harbourne and Giovanna Ilardi and Juan Migliore and Tomasz Szemberg and Justyna Szpond},
      year={2023},
      eprint={2312.04644},
      archivePrefix={arXiv},
      primaryClass={math.AG},
      url={https://arxiv.org/abs/2312.04644}, 
}

@article{Kettinger,
  title={The geproci property in positive characteristic},
  author={Kettinger, Jake},
  journal={Proceedings of the American Mathematical Society},
  volume={152},
  number={9},
  pages={3229-3242},
  year={2024},
  publisher={American Mathematical Society}
}

@article{Gorla,
title = {An algebraic approach for decoding spread codes},
author = {Elisa Gorla and Felice Manganiello and Joachim Rosenthal},
author = {Gorla, Elisa and Manganiello, Felice and Rosenthal, Joachim},
journal = {Advances in Mathematics of Communications},
volume = {6},
number = {4},
pages = {443-466},
year = {2012},
issn = {1930-5346},
doi = {10.3934/amc.2012.6.443},
url = {https://www.aimsciences.org/article/id/33546828-96a2-4bfd-87ef-90e7fe5d5082}
}

@article{BB,
  title={The construction of translation planes from projective spaces},
  author={Bruck, Richard and Bose, Raj},
  journal={Journal of Algebra},
  volume={1},
  number={},
  pages={85-102},
  year={1964},
  publisher={}
}

@article{M,
  title={Sets of Disjoint Lines in {$PG(3, q)$}},
  author={Mesner, Dale},
  journal={Canadian Journal of Mathematics},
  volume={19},
  number={},
  pages={273--280},
  year={1967},
  publisher={Cambridge University Press}
}

@Misc{M2,
          author = {Grayson, Daniel R. and Stillman, Michael E.},
          title = {Macaulay2, a software system for research in algebraic geometry},
          howpublished = {Available at \url{http://www2.macaulay2.com}}
        }

@Misc{groupnames,
          author = {Dokchitser, Tim},
          title = {GroupNames},
          howpublished = {Available at \url{https://people.maths.bris.ac.uk/~matyd/GroupNames/index.html}}
        }

@article{VlasovP,
    author = "Vlasov, Alexander Yurievich",
    title = "{Penrose Dodecahedron, Witting Configuration and Quantum Entanglement}",
    eprint = "2208.13644",
    archivePrefix = "arXiv",
    primaryClass = "quant-ph",
    doi = "10.12743/quanta.v13i1.276",
    journal = "Quanta",
    volume = "13",
    number = "1",
    pages = "38--46",
    year = "2024"
}

@article{Dolgachev,
author = {Dolgachev, Igor},
year = {2003},
month = {05},
pages = {},
title = {Abstract configurations in algebraic geometry},
journal = {The Fano Conference}
}

@manual{GAP,
    organization = {The GAP~Group},
    title        = {GAP -- Groups, Algorithms, and Programming, Version 4.14.0},
    year         = {2024},
    howpublished={\url{https://www.gap-system.org}}
    }

@book{Klein,
  title={Lectures on the icosahedron and the solution of equations of the fifth degree.},
  author={Klein, Felix},
  journal={},
  volume={},
  number={Translated into English by George Gavin Morrice.},
  pages={},
  year={1888},
  publisher={Tr{\"u}bner {\&} Co., Ludgate Hill}
}

@book{Brown,
  title={Topology and Groupoids: A Geometric Account of General Topology, Homotopy Types, and the Fundamental Groupoid},
  author={Brown, Ronald},
  journal={},
  volume={},
  number={},
  pages={},
  year={1988},
  publisher={Ellis Horwood Limited}
}

\end{document}